\documentclass{article}
\usepackage{amssymb, latexsym, amsmath, amsthm, graphicx}

\begin{document}

\newtheorem{thm}{Theorem}[section]
\newtheorem{prop}[thm]{Proposition}
\newtheorem{lem}[thm]{Lemma}
\newtheorem{cor}[thm]{Corollary}
\newtheorem{introthm}{Theorem}
\renewcommand{\theintrothm}{\Alph{introthm}}

\newtheorem{defn}[thm]{Definition}

\newcommand{\perisubs}{\{H_{i}\}_{i\in I}}
\newcommand{\cone}{\Gamma(G,X\cup\mathcal{H})}
\newcommand{\Cone}{\widehat{\Gamma}}
\newcommand{\unconed}{\Gamma(G,X)}
\newcommand{\coset}{gH_i}
\newcommand{\rtg}{\tau^{rel}(g)}
\newcommand{\parasub}{\gamma^{-1}H_{i}\gamma}
\newcommand{\hnnnf}{g_{0}t^{\varepsilon_{1}}g_{1}\ldots t^{\varepsilon_{n}}g_{n}}
\newcommand{\atbar}{\bar{a_{2}}}
\newcommand{\atobar}{\tau_{1}^{-1}a_{2}\tau_{1}}
\newcommand{\btobar}{\tau_{2}^{-1}b_{2}\tau_{2}}
\newcommand{\btbar}{\bar{b_{2}}}
\newcommand{\aobar}{\bar{a_{1}}}
\newcommand{\bobar}{\bar{b_{1}}}
\newcommand{\nn}{\|_{X}}
\newcommand{\mm}{\|_{\Sigma_{n}}}
\newcommand{\nf}{\|_{NF}}
\newcommand{\arr}{G}
\newcommand{\sizz}{\mathcal{L}}
\newcommand{\goo}{\gamma}
\newcommand{\goon}{\gamma_{0}}
\newcommand{\linlim}{Q}
\newcommand{\ylim}{(\linlim+4)\mathcal{L}}
\newcommand{\transversal}{T(A_{n})}
\newcommand{\NF}{\text{normal form}}
\newcommand{\CRNF}{\text{cyclically reduced normal form}}
\newcommand{\CRW}{\text{cyclically reduced word}}

\def\<{\langle}
\def\>{\rangle}

\title{Conjugacy Search Problem for Relatively Hyperbolic Groups}
\author{Zoe O'Connor\footnote{Heriot-Watt University, Edinburgh, U.K. Email: zao1@hw.ac.uk}}

\maketitle

\begin{abstract}
The asymptotic bound for a length-based attack on the Conjugacy Search Problem in relatively hyperbolic groups is cubic for hyperbolic elements and a ``small'' polynomial for parabolic elements, depending on the Conjugacy Search Problem for the peripheral subgroups. The bound for relatively hyperbolic groups in this paper is a significant improvement on previous work. 
\end{abstract}

\section{Introduction}
The \emph{Conjugacy Search Problem} (CSP) is the following: Given a group $G$, and two elements $a$ and $b$ which are conjugate in $G$, find an element $x\in G$ such that $x^{-1}ax=b$. Several key-agreement protocols (see, for example,~\cite{KAP}) are based on this problem. The aim of this paper is to put an upper bound $U$ on the minimum length of a conjugating element, so that the CSP in relatively hyperbolic groups can be solved by checking with all elements $x$ of length less than $U$ whether $x^{-1}axb^{-1}=1$ in $G$. The word problem in relatively hyperbolic groups was shown by Farb to have a ``fast'' solution:

\begin{thm}[\cite{FARB}, Theorem 3.7]
Suppose that a group $G$ is hyperbolic relative to a subgroup $H$, and H has word problem solvable in time $O(f(n))$. Then there is an algorithm which gives an $O(f(n)\log{n})$-time solution to the word problem in $G$.
\end{thm}

The \emph{Conjugacy Problem} is the following: Given two elements $a,b$ of a group $G$, determine whether $a$ is conjugate to $b$ in $G$. In~\cite{BUMA} Bumagin presented a proof that the Conjugacy Problem is solvable for relatively hyperbolic groups. Ji, Ogle and Ramsey used this paper to show that the CSP for relatively hyperbolic groups has a length bound which is a polynomial function of the lengths of the conjugate elements $a,b$~\cite{JOR}; a detailed study of this paper shows that this bound is a polynomial of degree $576n$, where $n$ is the degree of the polynomial bound for the CSP in the peripheral subgroups.

The results in this paper drastically improve this estimate, as follows:

\begin{introthm}
\label{thmA}
Let $G$ be a relatively hyperbolic group with generating set $X$, and suppose that the CSP in all peripheral subgroups can be bounded by a polynomial $\mathcal{P}$ of degree $n$. Let $a,b\in G$ be two elements which are conjugate in $G$, and let $\mathcal{L}=\max\{\|a\nn,\|b\nn\}$.
Then the CSP in $G$ is bounded by a polynomial function of $\mathcal{L}$ of degree $\max\{3, 2n+1\}$.
\end{introthm}

\section{Preliminaries}
There are several definitions of relatively hyperbolic groups in use; for the purposes of this paper we use the definition introduced by Farb in~\cite{FARB}. Let $G$ be a group with generating set $X$ and let $\mathcal{H}=\{H_{i}\}_{i\in I}$ be a set of subgroups in $G$ (called \emph{peripheral subgroups}). We call the Cayley graph $\Cone=\cone$ the \emph{coned-off Cayley graph} of $G$ with respect to $\mathcal{H}$. 

For a path $p$ in the Cayley graph $\Gamma=\unconed$ we denote the corresponding path in the coned-off graph $\Cone$ as $\hat{p}$. The path metrics in $\Gamma$ and $\Cone$ will be denoted by $d_{\Gamma}$ and $d_{\Cone}$ respectively. A geodesic, quasigeodesic, etc. in the coned-off graph will be called a \emph{relative} geodesic, quasigeodesic and so on. To be clear on the choice of generating set, the length of an element $x$ of $G$ with respect to the generating set $X$ will be denoted $\|x\nn$ and the relative length of $x$ with respect to $X\cup\mathcal{H}$ will be denoted $\|x\|_{X\cup\mathcal{H}}$. The length of a path $p$ in $\Gamma$ will be denoted $l_{\Gamma}(p)$ and the length of a path $\hat{p}$ in $\Cone$ will be denoted $l_{\Cone}(\hat{p})$. The label of a path $p$ in a Cayley graph will be written as $\phi(p)$ and will be identified with the element it represents in $G$. The centralizer of an element $g\in G$ will be written as $C_{G}(g)$. We denote the origin and terminus of a path $p$ by $p_{-}$ and $p_{+}$ respectively. 

The following vocabulary is borrowed from Osin's book~\cite{OSIN}; refer to this text for more background on relatively hyperbolic groups.

Two paths $p,q$ are called \emph{$k$-similar} if $d_{\Gamma}(p_{-},q_{-})\leq k$ and $d_{\Gamma}(p_{+},q_{+})\leq k$. We say that two paths $p,q$ are \emph{symmetric} if $\phi(p)\equiv\phi(q)$; i.e. if the two paths have identical labels. Given a pair of symmetric paths $(p,q)$ we call $g_{1}=(p_{-})^{-1}q_{-}$ and $g_{2}=(p_{+})^{-1}q_{+}$ the \emph{characteristic elements} of $(p,q)$. A symmetric pair of geodesics $(p,q)$ is said to be \emph{minimal} if for any other pair of symmetric geodesics $(p',q')$ with the same characteristic elements, the inequality $l_{\Cone}(\hat{p})\leq l_{\Cone}(\hat{p'})$ holds.

Let $(p,q)$ be a symmetric pair of paths. We say that the vertices $v_{1}$ of $p$ and $v_{2}$ of $q$ are \emph{synchronous vertices} if the path segments $[p_{-},v_{1}]$ and $[q_{-},v_{2}]$ have the same length. A subpath is called an \emph{$H_{i}$-component} if it is labelled by an element of $H_i$, and it is maximal in that respect (it is not contained in a larger subpath which is labelled by an element in $H_i$). 

Any vertex of a path $p$ which ``disappears'' in the coned-off graph $\Cone$ (that is, any vertex which is part of some $H_{i}$-component $s$ but is not equal to $s_-$ or $s_+$) is called \emph{non-phase}; all other vertices are called \emph{phase} vertices. 

Two $H_{i}$-components $s$ of $p$ and $t$ of $q$ are called \emph{synchronous} components if $s_-$ and $t_-$ are synchronous vertices. Two $H_{i}$-components $s$ of $p$ and $t$ of $q$ are called \emph{connected} components if there is a path in $\Gamma$ from $s_-$ to $t_-$ which is labelled by an element of $H_i$. When we speak of a single path $p$, we say that an $H_{i}$-component $s$ is \emph{isolated} if no distinct $H_{i}$ component of $p$ is connected to $s$ by a path in $H_{i}$. A path $p$ is called a \emph{path without backtracking} if every $H_{i}$-component of $p$ is isolated.

A finitely generated group $G$ is called \emph{weakly relatively hyperbolic} with respect to the subgroups $\mathcal{H}$ if the coned-off Cayley graph of $G$ with respect to $\mathcal{H}$ is hyperbolic with respect to the word metric. A further property is required for such a group to he called relatively hyperbolic:

\begin{defn}{Bounded Coset Penetration Property.}
Let $G$ be a weakly hyperbolic group relative to the subgroups $\perisubs$. Then $G$ is said to \emph{satisfy the Bounded Coset Penetration property (BCP)} if for any $\lambda$ there exists a constant $c(\lambda)$ such that the following conditions hold. Let $p,q$ be two relative $(\lambda,0)$-quasi-geodesics without backtracking, with the same endpoints;

\begin{enumerate}
\item If both $p$ and $q$ penetrate the same coset then they enter (and leave) the coset a distance at most $c(\lambda)$ apart.
\item If $p$ penetrates a coset $\coset$ which $q$ does not penetrate, then $p$ travels a distance at most $c(\lambda)$ in $\coset$.
\end{enumerate}

\end{defn}

\begin{defn}
A finitely generated group $G$ is said to be \emph{hyperbolic relative to its subgroups $\mathcal{H}$} (or simply \emph{relatively hyperbolic}) if it is weakly relatively hyperbolic with respect to $\mathcal{H}$ and it satisfies the Bounded Coset Penetration property.
\end{defn}

It has been shown~\cite{OSIN} that each peripheral subgroup $H_{i}$ is finitely generated, and that the set $\mathcal{H}=\perisubs$ is finite.

The next two results will be useful in the proof of Theorem~\ref{thmA}, particularly in establishing the size of certain constants. The first, called the \emph{Fellow Traveller Property}, is well known.

\begin{lem}
\label{FTP}
Let $(X,d)$ be a $\delta$-hyperbolic geodesic metric space, and let $\gamma_{1}:[0,T_{1}]\mapsto X$ and $\gamma_{2}:[0,T_{2}]\mapsto X$ be two $k$-similar geodesics. Then for any $t\leq max\{T_{1},T_{2}\}$ the points $\gamma_{1}(t)$ and $\gamma_{2}(t)$ are $(4\delta+3k)$-close.
\end{lem}

\begin{prop}
\label{o323q}
There is a polynomial $\varepsilon=\varepsilon(\lambda,c,k)$ such that for any two $k$-similar $(\lambda,c)$-quasi-geodesics without backtracking $p,q$ in $\cone$, the following conditions hold:

\begin{enumerate}
\item The sets of phase vertices of $p$ and $q$ are contained in the closed \\$\varepsilon$-neighbourhoods of each other.
\item Suppose that $s$ is an $H_{i}$-component of $p$ such that $d_{\Gamma}(s_{-},s_{+})>\varepsilon$, then there exists an $H_{i}$-component $t$ of $q$ which is connected to $s$.
\item Suppose that $s$ and $t$ are connected $H_{i}$-components of $p$ and $q$ respectively; then $max\{d_{\Gamma}(s_{-},t_{-}),d_{\Gamma}(s_{+},t_{+})\}\leq\varepsilon$.
\end{enumerate}

\end{prop}

This is Theorem 3.23 of~\cite{OSIN}; details of the proof can be found in that book. If $p$ and $q$ are geodesics then $\lambda=1$ and $c=0$; by analysing the details of~\cite{JOR} we see that $\varepsilon$ is in fact a quadratic function of $k$, and we can write it as $\varepsilon(k)$.

We state without proof a useful lemma from~\cite{OSIN}:

\begin{lem}
\label{o340}
Let $(\hat{p},\hat{q})$ be a minimal pair of symmetric relative geodesics in the Cayley graph $\cone$.
\begin{enumerate}
\item Suppose that, for some $i$, two $H_{i}$-components $s$ and $t$ of $\hat{p}$ and $\hat{q}$ respectively are connected. Then $s$ and $t$ are synchronous.
\item Let $u_{1},v_{1}$ and $u_{2},v_{2}$ be two pairs of synchronous vertices of $\hat{p}$ and $\hat{q}$ respectively. Then $(u_{1})^{-1}v_{1}\neq (u_{2})^{-1}v_{2}$.
\end{enumerate}
\end{lem}

\section{Conjugacy Search Problem for Relatively Hyperbolic Groups}
In this section $G$ will denote a group which is relatively hyperbolic with respect to the peripheral subgroups $\mathcal{H}=\perisubs$, with finite generating set $X$. Given two conjugate elements $a,b\in G$, our goal is to find an element $x$ which satisfies the equation $x^{-1}axb^{-1}=1$. Geometrically, we want to find a closed path $\Theta=\theta_{q}^{-1}\theta_{a}\theta_{p}\theta_{b}^{-1}$ in the Cayley graph of $G$ such that $\phi(\theta_{a})=a$, $\phi(\theta_{b})=b$ and $\phi(\theta_{p})=\phi(\theta_{q})=x$. We may assume that the path $\theta_{a}$ starts at the vertex labelled by the identity element. The subpaths $\theta_{p}$ and $\theta_{q}$ are symmetric and $\mathcal{L}$-similar, where $\mathcal{L}=\max\{\|a\nn,\|b\nn\}$. We want to find an upper bound on the length of the element $x$, so we will assume that $(\theta_{p},\theta_{q})$ is a minimal pair of symmetric geodesics, and we attempt to establish an upper bound on the $\Gamma$-length of these geodesics.

\begin{figure}[ht]
\centering\includegraphics[scale=.45]{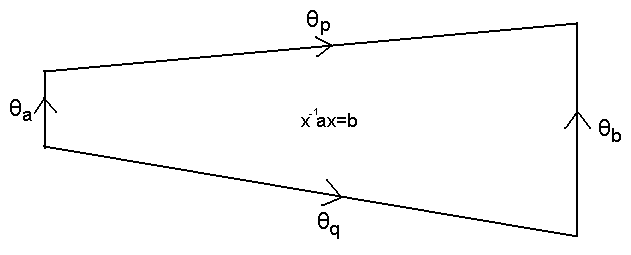}
\caption{Conjugacy diagram showing $\Theta$ in $\unconed$}
\label{F211:conj}
\end{figure}

\begin{defn}
Let $G$ be a group which is hyperbolic relative to $\mathcal{H}=\perisubs$. An element $g\in G$ is \emph{parabolic} if it is conjugate to some element of one of the peripheral subgroups $\mathcal{H}$, otherwise it is called \emph{hyperbolic}.
\end{defn}

\begin{lem}\label{tor1}
Let $a\in G$ be conjugate to an element $b$ by some conjugating element $x$ of minimal length. If $(u,v)$ is a pair of synchronous vertices on $(\theta_{p},\theta_{q})$ with $d_{\Cone}((\theta_{p})_{\pm},u)>\mathcal{L}+2\delta$, then $d_{\Cone}(u,v)\leq 4\delta$.
\begin{proof}
As usual we assume that $\theta_{p}$ and $\theta_{q}$ are chosen to be minimal. We parametrize $\theta_{p}$ and $\theta_{q}$ so that $\theta_{p}(i)$ is the $i^{th}$ vertex along the path $\theta_{p}$; likewise with $\theta_{q}$. Let $u=\theta_{p}(t)$; then $v=\theta_{q}(t)$.

By the $2\delta$-thinness of quadrilaterals in hyperbolic spaces, there is some vertex $\theta_{p}(t')$ which is $2\delta$-close to $\theta_{q}(t)$. Suppose without loss of generality that $t'\geq t$. Then $\theta_{p}(t)$ is $2\delta$-close to $\theta_{r}(t')$, where $\theta_{r}=a\theta_{p}$ (see Fig.~\ref{F233:conjug}).

\begin{figure}[ht]
\centering\includegraphics[scale=.5]{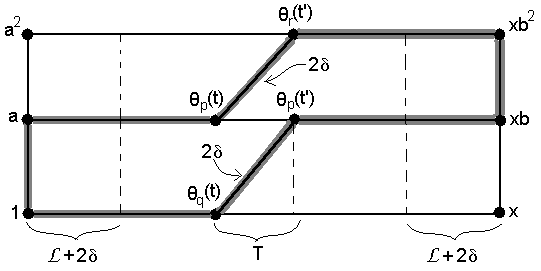}
\caption{Conjugacy diagram for Lemma \ref{tor1}}
\label{F233:conjug}
\end{figure}

If $T=t'-t>2\delta$ then the paths \[[a,\theta_{p}(t)][\theta_{p}(t),\theta_{r}(t')][\theta_{r}(t'),xb^{2}]\] and \[[1,\theta_{q}(t)][\theta_{q}(t),\theta_{p}(t')][\theta_{p}(t'),xb]\] are shorter than $\theta_{p}$ and $\theta_{q}$, and they conjugate $a$ and $b$ (as highlighted in grey on the diagram), which contradicts our assumption that $(\theta_{p},\theta_{q})$ is a minimal pair of synchronous geodesics. It now follows from the triangle inequality that $d_{\Cone}(u,v)\leq T+2\delta\leq 4\delta$.
\end{proof}
\end{lem}

This shows that there is a ``middle'' section of $\theta_{p}$ with $\Cone$-length bounded by the number of distinct words in $X\cup\mathcal{H}$ of length $4\delta$. We make use of the following:

\begin{lem}[\cite{OSIN}, Lemma 3.41] \label{o341}
Let $(p,q)$ be a minimal pair of symmetric geodesics in $\unconed$ such that
\[\max\{d_{\Cone}(\hat{p}_{-},\hat{q}_{-}),d_{\Cone}(\hat{p}_{+},\hat{q}_{+})\}\leq k\]
and let $v_{1}, v_{2}$ be synchronous vertices on $p$ and $q$ respectively such that
\[\min\{d_{\Cone}(\hat{p}_{-},v_{1}),d_{\Cone}(\hat{p}_{+},v_{1})\}\geq 2E,\]
where $E=4\delta+3k$ is the constant from Lemma~\ref{FTP}. Then $$d_{\Gamma}(v_{1},v_{2})\leq 6MLE^{2},$$
where $ML>1$ is a constant (see~\cite{OSIN}, Convention 3.1 for details).
\end{lem}

Setting $k=4\delta$ and combining this with the previous lemmas we have the following:

\begin{lem}\label{torl}
Let $a\in G$ be conjugate to an element $b\in G$. Then there exists $x\in G$ such that $a=x^{-1}bx$ and $$\|x\|_{X\cup\mathcal{H}}\leq 2(\mathcal{L}+34\delta)+|X|^{6ML(16\delta)^{2}},$$
where $\mathcal{L}=\max\{\|a\nn,\|b\nn\}$.
\begin{proof}
We have established in Lemma~\ref{tor1} that synchronous vertices which are a $\Cone$-distance of at least $\mathcal{L}+2\delta$ from either end of $\hat{\theta_{p}}$ and $\hat{\theta_{q}}$ respectively are a $\Cone$-distance of at most $4\delta$ apart. Let us call these middle sections $\hat{\theta_{p}'}$ and $\hat{\theta_{q}'}$. Then we can use Lemma~\ref{o341} with $k=4\delta$ to prove that if $(u,v)$ is a pair of synchronous vertices on the paths $\theta_{p}$ and $\theta_{q}$ such that $d_{\Cone}(u,(\theta_{p}')_{\pm})\geq 2E=2(4\delta+3(4\delta))=32\delta$ then $d_{\Gamma}(u,v)\leq 6ML(16\delta)^{2}$. Thus the length of the section of $\theta_{p}$ which is a $\Cone$-distance of $\mathcal{L}+34\delta$ from either end of $\theta_{p}$ has $\Gamma$-length at most $|X|^{6ML(16\delta)^{2}}$ by the argument that if there are two pairs of synchronous vertices which are joined by a geodesic of the same label, then we can shorten the path $\Theta$ by ``cutting out'' the section between these two geodesic paths and joining the remaining parts together along these geodesics.

We conclude that $\hat{\theta_{p}}$ has a $\Cone$-length of $l_{\Cone}(\hat{\theta_{p}})\leq 2(\mathcal{L}+34\delta)+|X|^{6ML(16\delta)^{2}}$, as illustrated in Fig~\ref{F344:length}, in which lengths are $\Cone$-lengths unless otherwise stated.
\begin{figure}[ht]
\centering\includegraphics[scale=.35]{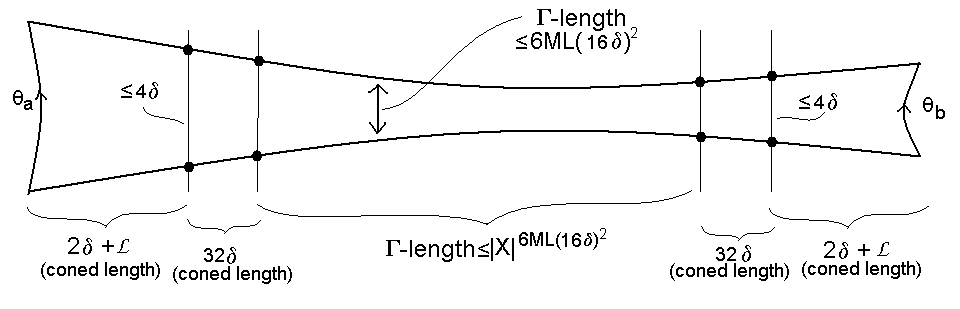}
\caption{Length diagram for Lemma \ref{torl}}
\label{F344:length}
\end{figure}
\end{proof}
\end{lem}

The following is drawn from results in~\cite{OSIN}:

\begin{lem}\label{hyptravel}
Let $a$, $b$ be conjugate hyperbolic elements of $G$, with a conjugating element $x$ of minimal length. The $\Gamma$-distance through which the associated paths $\theta_{p}$ and $\theta_{q}$ of the cycle $\Theta$ travel in each $H_{i}$-coset is bounded above by the quadratic function $\varepsilon(\mathcal{L})$, where $\mathcal{L}=\max\{\|a\nn,\|b\nn\}$, and $\varepsilon$ is the quadratic from Proposition~\ref{o323q}.
\begin{proof}
Consider the closed cycle $\Theta$ in $\unconed$. Proposition~\ref{o323q} states that if $s$ is an $H_{i}$-component of $\theta_{p}$ with $l_{\Gamma}(s)>\varepsilon$, then there exists an $H_{i}$-component $t$ of $\theta_{q}$ which is connected to $s$ by a path labelled by $h\in H_{i}$. Furthermore, by Lemma~\ref{o340}, since $\theta_{p}$ and $\theta_{q}$ are minimal and symmetric, and $s$ and $t$ are connected, then these two components are synchronous. Consequently $a$ and $b$ are conjugate to $h\in H_{i}$, but we are assuming that $a$ and $b$ are hyperbolic; a contradiction. Hence $l_{\Gamma}(s)\leq\varepsilon$.
\end{proof}
\end{lem}

\begin{thm}
\label{rhghyper}
Let $a$ and $b$ be conjugate hyperbolic elements of the relatively hyperbolic group $G$. Then there exists $x\in G$ such that $x^{-1}ax=b$ and $\|x\nn$ is bounded above by a cubic polynomial in $\mathcal{L}$.
\begin{proof}
Lemma~\ref{hyptravel} shows that $\theta_{p}$ travels a $\Gamma$-distance of no more than $\varepsilon$ in each coset it penetrates. By Lemma~\ref{torl} we know that there is a `middle section' of $\theta_{p}$ which has $\Gamma$-length bounded by the constant $|X|^{6ML(16\delta)^{2}}$. Either side of this section is a subpath of $\theta_{p}$ which has $\Cone$-length bounded by $34\delta+\mathcal{L}$. Hence
\[\|x\nn=l_{\Gamma}(\theta_{p})\leq 2(34\delta+\mathcal{L})\varepsilon(\mathcal{L})+|X|^{6ML(16\delta)^{2}}\]
which is a cubic polynomial in $\mathcal{L}$.
\end{proof}
\end{thm}

\begin{lem}
\label{pring}
Let $a$ and $b$ be conjugate parabolic elements in $G$ with respect to $\perisubs$, and suppose that the conjugacy search problem in each of the subgroups $H_{i}$ is bounded above by a polynomial $\mathcal{P}$ of degree $n$. Then the paths $\theta_{p}$ and $\theta_{q}$ each travel a polynomially bounded distance, of degree $2n$, in each coset they penetrate.
\begin{proof}
Choose a peripheral subgroup $H_{i}$ and consider the set $\Sigma=\{(u_{j},v_{j}):j=1,\ldots,n\}$ of all synchronous phase vertices on $\theta_{p}$ and $\theta_{q}$ respectively which are each joined by a geodesic path in $\cone$ labelled by $h_{j}\in H_{i}$, such that $\theta_{p}$ reaches each $u_{j}$ in ascending order. Divide the quadrilateral $\Theta$ into $n+1$ ``cells'' using these paths $\{h_{j}\}_{j=1}^{n}$. For notational ease, let $h_{0}:=\theta_{a}$ and $h_{n+1}:=\theta_{b}$.

\begin{figure}[ht]
\centering\includegraphics[scale=.5]{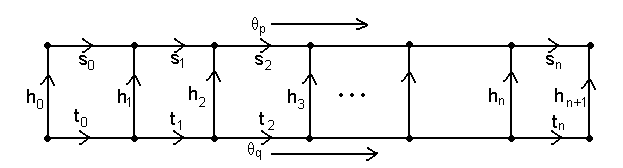}
\caption{Illustration of the ``cells'' in Lemma \ref{pring}}
\label{F453:cells}
\end{figure}

The segment $s_{j}:=[(h_{j})_{+},(h_{j+1})_{+}]$ of $\theta_{p}$ in each cell will fall into one of two categories. In the first case it is an $H_{i}$-component, in which case so is $t_{j}:=[(h_{j})_{-},(h_{j+1})_{-}]$. Then $h_{j}$ and $h_{j+1}$ are conjugate in $H_{i}$, and the $\Gamma$-length of $s_{j}$ and $t_{j}$ will be bounded by $\mathcal{P}(\max\{l_{\Gamma}(h_{j}),l_{\Gamma}(h_{j+1})\})\leq\mathcal{P}(\varepsilon)$, where $\varepsilon$ is the quadratic polynomial from Proposition~\ref{o323q}. Note that if $s_{0}$ is an $H_{i}$-component then $a=h^{-1}h_{1}h$ for some $h\in H_{i}$ and hence $u\in H_{i}$; likewise if $s_{n}$ is an $H_{i}$-component then $b\in H_{i}$.

The second case is that $s_{j}$ is not an $H_{i}$-component of $\theta_{p}$, although it may contain an $H_{i}$-component $h'$ of $\theta_{p}$ which, by our choice of the paths $h_{j}$, will not be connected to the synchronous $H_{i}$-component of $\theta_{q}$. Since $(s_{j},t_{j})$ is a minimal pair of synchronous geodesics with characteristic elements $h_{j}$ and $h_{j+1}$, we can use Lemma~\ref{o340} to see that if $h'$ is connected to any $H_{i}$-component of $\theta_{q}$ then these two components must be synchronous. Then by Proposition~\ref{o323q}, as $h'$ is an isolated component, its $\Gamma$-length is bounded by the quadratic $\varepsilon(\mathcal{L})$.

We conclude that for parabolic $a,b$ the $\Gamma$-length of any $H_{i}$-component of $\theta_{p}$ is bounded by a polynomial  $\mathcal{M}(\mathcal{L})=\max\{\varepsilon(\mathcal{L}),\mathcal{P}(\varepsilon(\mathcal{L}))\}$.
\end{proof}
\end{lem}

\begin{thm}
\label{rhgpara}
Let $a$ and $b$ be conjugate parabolic elements of the relatively hyperbolic group $G$. Suppose that the CSP in each peripheral subgroup can be bounded by a polynomial $\mathcal{P}(\mathcal{L})$ of degree $n$. Then there exists an element $x\in G$ such that $a^{x}=b$ and the $\Gamma$-length of $x$ is bounded by a polynomial of degree $2n+1$.
\begin{proof}
As in the hyperbolic case, the $\Gamma$-length of $\theta_{p}$ is bounded by $$2(34\delta+k)\mathcal{M}+|X|^{6ML(16\delta)^{2}}$$
where $\mathcal{M}$ is the polynomial from Lemma~\ref{pring}. Then $l_{\Gamma}(\theta_{p})$ is bounded above by a polynomial in $\mathcal{L}$ of degree $2n+1$.
\end{proof}
\end{thm}

\setcounter{introthm}{0}
\begin{introthm}
Let $G$ be a relatively hyperbolic group with generating set $X$, and suppose that the CSP in all peripheral subgroups can be bounded by a polynomial $\mathcal{P}$ of degree $n$. Let $a,b\in G$ be two elements which are conjugate in $G$, and let $\mathcal{L}=\max\{\|a\nn,\|b\nn\}$.
Then the CSP in $G$ is bounded by a polynomial function of $\mathcal{L}$ of degree $\max\{3, 2n+1\}$.
\begin{proof}
The result follows from Theorem~\ref{rhghyper} and Theorem~\ref{rhgpara}.
\end{proof}
\end{introthm}

\begin{cor}
If $G$ is a limit group then for any pair of conjugate elements $a,b$ we can find a conjugating element $x$ of length at most $\mathcal{P}$, where $\mathcal{P}$ is a cubic polynomial in $\mathcal{L}=\max\{\|a\nn,\|b\nn\}$.
\begin{proof}
Limit groups are hyperbolic relative to their maximal non-cyclic abelian subgroups~\cite{FD}, and the conjugacy search problem in abelian groups is trivial. Hence $n=0$ and the cubic bound for the hyperbolic case gives the asymptotic upper bound for $\|x\nn$.
\end{proof}
\end{cor}

In fact, this is not a best-possible bound for limit groups. A result currently in preparation by the author of this paper will show that the CSP for limit groups has an asymptotically linear upper bound, as does the multiple CSP for limit groups.
\section{Open Problems}
This paper does not provide a time bound on solving the CSP for relatively hyperbolic groups. The emphasis is on finding a length bound; so far an algorithm has not been considered.

The conjugacy search problem can be generalized as follows: Given two lists of elements, $A=[a_{1},\ldots,a_{m}]$ and $B=[b_{1},\ldots,b_{m}]$, which are conjugate in a group $G$, find an element $x$ such that $x^{-1}a_{i}x=b_{i}$ for all $i=1,\ldots,m$. This is known as the \emph{multiple conjugacy search problem}. What is the asymptotic bound for a length-based attack on the multiple CSP for relatively hyperbolic groups? 

Bridson and Howie~\cite{BHo2005} showed that the multiple CSP for hyperbolic groups has a linear asymptotic bound, and their argument is based on the Cayley graph of the hyperbolic groups. The problem with using a similar argument for relatively hyperbolic groups is that the argument relies on putting a finite upper bound on the size of the centralizers by looking at the action of certain elements on the Cayley graph to see if the centralizers intersect at a finite number of points. However if we look at the action of a relatively hyperbolic group on $\Cone$ then a ball of bounded radius in $\Cone$ contains an infinite number of vertices.
\section*{Acknowledgements}
This work was funded by the EPSRC DTA grant EP/P504945/1.

\bibliographystyle{plain}
\bibliography{Mainbibliography}
\end{document}